\documentclass[12pt]{amsart}

\usepackage{a4wide}
\usepackage{amsmath}
\usepackage{amssymb}

\newtheorem{stw}{{\bf Proposition}}

\newenvironment{dow}{\noindent{\bf\em Proof.}\,}
{\hspace*{\fill}\(\boxtimes\)\newline}

\newcounter{liczprz}
\newenvironment{prz}{\refstepcounter{liczprz}{\noindent
\bf Example~\theliczprz.}}{\hspace*{\fill}\(\lozenge\)\newline}

\newcommand{\B}{\mathcal{B}}
\newcommand{\C}{\mathcal{C}}
\newcommand{\Pow}{\mathcal{P}}

\newcommand{\sprz}[1]{\overline{#1}}
\newcommand{\bdr}[1]{\partial({#1})}
\newcommand{\bd}{\partial}
\newcommand{\op}{\operatorname{op}}

\begin{document}

\title[Axiomatization]
{Axiomatization of topological space in terms of the operation of
boundary}

\author[Le\'{s}niak]{K. Le\'{s}niak}
\date{May 9, 2006}



\address{Faculty of Mathematics and Computer Science, Nicolaus Copernicus
University, ul. Chopina 12{\slash}18, 87-100 Toru\'{n}, Poland}

\email{much@mat.uni.torun.pl}

\begin{abstract}
We present the set of axioms for topological space with the
operation of boundary as primitive notion.
\end{abstract}

\maketitle%

Let \(X\) denote space (without any ascribed structure) and
\(\Pow(X)\) family of its subsets. We say that \(\sprz{(\cdot)}:
\Pow(X) \to \Pow(X)\) is {\it closure operation} if for any sets
\(A,B\subset X\)
\begin{enumerate}
 \item[(\(\delta\)-1)] \(\sprz{\emptyset}=\emptyset\),
 \item[(\(\delta\)-2)] \(\sprz{\sprz{A}} \subset \sprz{A}\),
 \item[(\(\delta\)-3)] \(\sprz{A\cup B}\subset
 \sprz{A}\cup\sprz{B}\),
 \item[(\(\delta\)-4)] \(A\subset B \Rightarrow
 \sprz{A}\subset\sprz{B}\),
 \item[(\(\delta\)-5)] \(A\subset\sprz{A}\).
\end{enumerate}
 We say that \(\partial:\Pow(X) \to \Pow(X)\) is {\it operation of boundary} if
for any sets \(A,B\subset X\)
\begin{enumerate}
 \item[(\(\beta\)-1)] \(\bd{\emptyset}=\emptyset\),
 \item[(\(\beta\)-2)] \(\bd{\bd{A}} \subset \bd{A}\),
 \item[(\(\beta\)-3)] \(\bdr{A\cup B}\subset
 \bd{A}\cup\bd{B}\),
 \item[(\(\beta\)-4)] \(A\subset B \Rightarrow
 \bd{A}\subset B\cup\bd{B}\),
 \item[(\(\beta\)-5)] \(\bd{A}=\bdr{X\setminus A}\).
\end{enumerate}
Axiom (\(\beta\)-5) can be still weakened to
\begin{enumerate}
 \item[(\(\beta\)-5')] \(\bd{A}\subset\bdr{X\setminus A}\).
\end{enumerate}

Below we give natural correspondence between these notions. Define
\(\Phi:{\Pow(X)}^{\Pow(X)}\to{\Pow(X)}^{\Pow(X)}\),
\(\forall_{\op\in{\Pow(X)}^{\Pow(X)}} \;\forall_{A\in\Pow(X)}
\;\;[\Phi(\op)]\,(A) \doteq \op(A)\cap\op(X\setminus A)\), and
\(\Psi:{\Pow(X)}^{\Pow(X)}\to{\Pow(X)}^{\Pow(X)}\),
\(\forall_{\op\in{\Pow(X)}^{\Pow(X)}} \;\forall_{A\in\Pow(X)}
\;\;[\Psi(\op)]\,(A) \doteq A\cup \op(A)\).

\begin{stw}[boundary via closure]
If \(\sprz{(\cdot)}\) is closure operation, then
\(\Phi(\sprz{(\cdot}))\) is operation of boundary.
\end{stw}
\begin{dow}
All calculations are standard so we show for example that
\(\partial\doteq\Phi(\sprz{(\cdot)})\) satisfies (\(\beta\)-4). If
\(A\subset B\), then \(\sprz{A}\subset\sprz{B}\) in view of
(\(\delta\)-4). Hence \(\bd{A}=\sprz{A}\cap\sprz{X\setminus A}
\subset\sprz{A}\subset\sprz{B}\). Further
\[\sprz{B}\subset
\sprz{B}\cup(X\setminus\sprz{X\setminus B}) =
(\sprz{B}\cap\sprz{X\setminus B}) \cup (X\setminus\sprz{X\setminus
B}) \stackrel{(\ast)}{\subset} (\sprz{B}\cap\sprz{X\setminus B})
\cup B = B\cup \bd{B},\] where inclusion \((\ast)\) is due to
(\(\delta\)-5).
\end{dow}

\begin{stw}[closure via boundary]\label{domknieciepoprzezbrzeg}
If \(\partial\) is operation of boundary, then \(\Psi(\partial)\)
is closure operation.
\end{stw}
\begin{dow}
Since most calculations are straightforward we only demonstrate
that \(\sprz{(\cdot)}\doteq\Psi(\partial)\) fulfills
(\(\delta\)-2) and (\(\delta\)-4).

ad (\(\delta\)-2): \(\sprz{\sprz{A}} = \sprz{A\cup\bd{A}} =
(A\cup\bd{A}) \cup \bdr{A\cup\bd{A}} \stackrel{(\ast)}{\subset}
A\cup\bd{A}\cup\bd{\bd{A}} \stackrel{(\ast\ast)}{\subset}
A\cup\bd{A}=\sprz{A}\), where \((\ast)\) uses (\(\beta\)-3) and
\((\ast\ast)\) uses (\(\beta\)-2).

ad (\(\delta\)-4): if \(A\subset B\), then \(\bd{A}\subset
B\cup\bd{B}\) by (\(\beta\)-4). Hence
\(\sprz{A}=A\cup\bd{A}\subset B\cup\bd{B}=\sprz{B}\).
\end{dow}

Denote by \(\C,\B\subset{\Pow(X)}^{\Pow(X)}\) the family of
closure and respectively boundary operations.

\begin{stw}[equivalence of definitions]\label{odpowiedniosc}
The correspondences \(\Phi:\C\to\B\) and \(\Psi:\B\to\C\) are
mutually inverse. In particular \(\Phi\) and \(\Psi\) are
bijections.
\end{stw}
\begin{dow}
Let \(\sprz{(\cdot)}\in \C\), \(A\in\Pow(X)\) and
\(\bd\doteq\Phi(\sprz{(\cdot)})\). Then
\(\Psi(\Phi(\sprz{(\cdot)}))\,(A) = A\cup \bd{A} = A\cup
(\sprz{A}\cap\sprz{X\setminus A}) \stackrel{(\ast)}{\subset}
\sprz{A}\), where \((\ast)\) uses (\(\delta\)-5). To get the
reverse inclusion in \((\ast)\) axiom (\(\delta\)-5) is used
again: \(\sprz{A}\setminus A\subset X\setminus A\subset
\sprz{X\setminus A}\), so \(\sprz{A} = A\cup(\sprz{A}\setminus A)
\subset A\cup \sprz{X\setminus A}\).

Now let \(\partial\in \B\), \(A\in\Pow(X)\) and
\(\sprz{(\cdot)}\doteq\Psi(\bd)\). Then \(\Phi(\Psi(\bd))\,(A) =
\sprz{A}\cap \sprz{X\setminus A} = (A\cup\bd{A}) \cap
\left((X\setminus A)\cup \bdr{X\setminus A}\right)
\stackrel{(\ast)}{=} (A\cup\bd{A}) \cap \left((X\setminus A)\cup
\bd{A}\right) = (A\cap(X\setminus A))\cup\bd{A} = \bd{A}\), where
\((\ast)\) uses (\(\beta\)-5).
\end{dow}

Observe that axiom (\(\beta\)-5) used to prove
Proposition~\ref{odpowiedniosc} is not exploited in the proof of
Proposition~\ref{domknieciepoprzezbrzeg}. Nevertheless this axiom
is indispensable as claimed by

\begin{stw}[logical independence]
The system of axioms (\(\beta\)-1) -- (\(\beta\)-5) is logically
independent.
\end{stw}

We split the verification of the above proposition in the series
of examples.
\newline

\begin{prz}
 Put \({\bd}_1(A)\doteq X\) for any \(A\subset X\). Then
\({\bd}_1\) fulfills all axioms of boundary except (\(\beta\)-1).
\end{prz}

\begin{prz}
 Let \(X\doteq \mathbb{N}\) and
 \[{\bd}_2(A)\doteq \left\{x\in\mathbb{N}\,:\,\inf_{a\in A} |x-a| =1 \vee
 \inf_{b\in {\mathbb{N}\setminus A}} |x-b| =1\right\}\]
 for any \(A\subset X\). Then
\({\bd}_2\) fulfills all axioms of boundary except (\(\beta\)-2).
\end{prz}

\begin{prz}
 Let \(X\doteq\{1,2,3\}\) and \[{\bd}_3(A)\doteq
 \left\{\begin{array}{ll}
 \emptyset, & \;\mbox{if}\; A=\emptyset \;\mbox{or}\; X,\\
 A, & \;\mbox{if}\; A= \{1\} \;\mbox{or}\; \{2\} \;\mbox{or}\; \{3\},\\
 X\setminus A, & \;\mbox{if}\; A= \{1,2\} \;\mbox{or}\; \{2,3\}
 \;\mbox{or}\; \{1,2\},\\
 \end{array}\right.\]
for any \(A\subset X\). Then \({\bd}_3\) fulfills all axioms of
boundary except (\(\beta\)-3).
\end{prz}

\begin{prz}
 Let \(X\doteq\{1,2,3\}\) and \[{\bd}_4(A)\doteq
 \left\{\begin{array}{ll}
 \emptyset, & \;\mbox{if}\; A=\emptyset \;\mbox{or}\; X,\\
 \{2\}, & \;\mbox{if}\; A= \{1\} \;\mbox{or}\; \{2,3\},\\
 \{1\}, & \;\mbox{if}\; A= \{2\} \;\mbox{or}\; \{1,3\},\\
 \{1,2\}, & \;\mbox{if}\; A= \{3\} \;\mbox{or}\; \{1,2\}.\\
 \end{array}\right.\]
for any \(A\subset X\). Then \({\bd}_4\) fulfills all axioms of
boundary except (\(\beta\)-4).
\end{prz}

\begin{prz}
 Fix \(x_0\in X\) and put \({\bd}_5(A)\doteq A\cup\{x_0\}\)
for every nonempty \(A\subset X\) and
\({\bd}_5(\emptyset)=\emptyset\). Then \({\bd}_5\) fulfills all
axioms of boundary except (\(\beta\)-5).
\end{prz}

In the context of the last example we remark that any closure
operation satisfies (\(\beta\)-1) -- (\(\beta\)-4) but never
(\(\beta\)-5).

\end{document}